\documentclass[12pt]{article}

\title{\bf{Restriction of stable bundles on a jacobian of genus $2$ to an embedded curve}} 

\author{\bf{Cristian Anghel}}  

\newtheorem{th1}{Theorem}[section]
\newtheorem{rem1}[th1]{Remark}
\newtheorem{th2}[th1]{Theorem}
\newtheorem{th3}[th1]{Theorem}
\newtheorem{rem2}[th1]{Remark}
\newtheorem{th6}[th1]{Theorem}
\newtheorem{th7}[th1]{Theorem}
\begin{document}

\maketitle

\bigskip
{This paper is dedicated to the memory of Professor Gheorghe Galbura.}
\bigskip

\bigskip
\small{The aim of this note is to describe the restriction map from the moduli space of
stable rank $2$ bundle with small $c_2$ on a jacobian $X$ of dimension $2$, to the moduli space 
of stable rank $2$ bundles on the corresponding genus $2$ curve $C$ embedded in $X$. }
\bigskip


\begin{center}
\section{Introduction}
\end{center}

Let $C$ a smooth curve of genus $2$ and $X$ his jacobian wich is a smooth projective algebraic surface.
We denote by $M_{(2, \ C, \ i)}$ for $i=1 \ or \ 2$ the moduli space of rank $2$ bundle on $X$ with $c_1=C$ 
and $c_2=i$. Also we denote by $M_{(2, \ K)}$ the moduli space of rank $2$ bundle on $C$ with determinant $K$ 
i.e. the canonical class of $C$. Obviously, for any $E \in M_{(2, \ C, \ i)}$ the restriction $E_{\mid C}$ is a 
rank $2$ bundle on $C$ with determinant $K$. 

\newpage

The natural questions wich appear are the followings: is $E_{\mid C}$
a stable (or at least semi-stable) bundle on $C$ and if yes, what is the induced map $M_{(2, \ C, \ i)}\longrightarrow M_{(2, \ K)}$ ? As we shall see, the answer depend on $i$: for $i=1$, the restriction is semi-stable, but for $i=2$ and $E$ generic 
in $M_{(2, \ C, \ 2)}$ the restriction is stable. Also, in the second case we can describe for some non-generic bundles $E$  
what is the restriction $E_{\mid C}$.

\begin{center}
\section{Previously known results}
\end{center}

For $X$ the jacobian of a genus $2$ curve $C$, we denote by $F_0={\cal{O}}(C)\otimes{\cal{J}}_{0}$, where ${\cal{J}}_{0}$ is 
the sheaf of ideals of the origin of $X$. Also, using $F_0$ we can construct a unique extension 
$0\longrightarrow {\cal{O}}_X\longrightarrow F_{-1}\longrightarrow F_0\longrightarrow 0$ wich has $c_1={\cal{O}}_X(C)$ and 
$c_2=1$. The first result we need is the following, proved in \cite{muk2}:

\begin{th1}
For any rank $2$ bundle $E$ on $X$ with $c_1={\cal{O}}_X(C)$ and 
$c_2=1$ there are uniques $x , \ y\in X$ such that $E\simeq {T_x}^*F_{-1}\otimes P_y$, where ${T_x}^*$ is the pull-back 
by the x-translation and $P_y$ is the line bundle on $X$ wich correspond to $y$ by the canonical isomorphism 
$X \longrightarrow \widehat X$ defined by the principal polarisation $C$. As consequence the moduli space is isomorphic 
with $X\times X$.
\end{th1}

It is very easy to verify that the condition for $E$ the have $det(E)={\cal O}_X(C)$ is that $x=-2y$; so we have the following: 

\begin{rem1}
The moduli space of rank $2$ bundles on $X$ with $c_1={\cal O}_X(C)$ and $c_2=1$ is isomorphic with $X$. 
\end{rem1}

For the moduli space on $C$ we need the following theorem proved in \cite{nr}:

\begin{th2}  \label{1}
Let $F$ a semi-stable rank $2$ bundle on $C$ with determinant equal with the canonical class of $C$, and $x_0$ a Weierstrass 
point of $C$. 
Let $D_F=\{\xi \in Pic^1(C) \ \mid H^{0}(\xi \otimes F\otimes {\cal O}(-x_0))\not=0\}$. With these notations, $D_F$ is a divisor 
of the linear system $\mid 2C\mid $ on $Pic^1(C)$ and the map $F\longrightarrow D_F$ is an isomorphism between the moduli 
space of rank two bundles with canonical determinant and $\mathbf{P}^3$.

\end{th2}  

For the case $c_2=2$ we need the following result proved in \cite{muk1} and \cite{u}:

\begin{th3}
 
$M_{(2, \ C, \ 2)}$ is isomorphic with $X\times Hilb^3(X)$, and for any $E\in M_{(2, \ C, \ 2)}$ there exist an unique 
exact sequence of the form:
\begin{center}
$0\longrightarrow {T_x}^*{\cal O}_{X}(-C)\longrightarrow H\longrightarrow E\longrightarrow 0$
\end{center}
where $H$ is an homogenous rank $3$ bundle on $X$.
\end{th3}

By \cite{muk2} a generic homogenous rank $3$ bundle has the form $P_a\oplus P_b\oplus P_c$ with $a\not=b\not=c$ and it is clear that the condition for $E$ the have $det(E)={\cal O}_X(C)$ is that $x=-a-b-c$; so we have the following: 

\begin{rem2}
The moduli space of rank $2$ bundles on $X$ with $c_1={\cal O}_X(C)$ and $c_2=2$ is birational with $Sym^3(X)$. 
\end{rem2}

\

\begin{center}
\section{The restriction theorems}
\end{center}

Using the previous notations we have the followings:

\begin{th6}
For generic $y\in X$ the restriction $E_{\mid C}$ of $E\simeq {T_{-2y}}^*F_{-1}\otimes P_y$ is semi-stable but not stable.
The rational restriction map $X--\rightarrow \mathbf{P}^3$ is the quotient by the natural involution of $X$ and the image 
is the Kummer surface.

\end{th6}

\begin{th7}
For generic $E \in Hilb^3(X)$ the restriction $E_{\mid C}$ is stable. The restriction $E_{\mid C}$, viewed in $\mathbf{P}^3=\mid 2C\mid $ is the unique divisor of $\mid 2C\mid $ wich contain the $3$ points $a, \ b, \ c$ of the corresponding $H$. Also, the fiber over a point $C'\in \mid 2C\mid $ is birational with $Hilb^3(C')$.
\end{th7}

The main idea in the proof of the previous theorems is to obtain an explicit description of $D_{E_{\mid C}}$ for generic $E$ in
the corresponding moduli space. In the first case for generic $y\in X$ and $E\simeq {T_{-2y}}^*F_{-1}\otimes P_y$ we obtain that 
$D_{E_{\mid C}}$ is the union of the two translate of $C$ by $y$ and $-y$. For $c_2=2$ and generic $E$, $D_{E_{\mid C}}$ is the hyperplane wich pass by the 3 points wich determine the homogenous bundle $H$ associated with $E$ by \ref{1} above.
The full details will appear elsewere.

\bigskip

\noindent 
Cristian Anghel \\
Department of Mathematics\\
Institute of Mathematics of the Romanian Academy\\
Calea Grivitei nr. 21 Bucuresti Romania\\
email:\textit{Cristian.Anghel@imar.ro}

\bigskip
\begin{thebibliography}{99}
\bibitem{muk1} S. Mukai:
               {\sl Fourier functor and its applications to the moduli of bundles on an abelian variety}, Algebraic Geometry, 
               Sendai, 1985.

\bibitem{muk2} S. Mukai:
               {\sl Duality between $D(X)$ and $D(\widehat X)$ with applications to Picard sheaves}, Nagoya Math. J. 81 1981 
               153-173.
\bibitem{nr}   Narasimhan, Ramanan:
               {\sl Moduli of vector bundles on a compact riemann surface}, Ann. of Math. 1969 14-51.
\bibitem{u}    H. Umemura:
               {\sl Moduli spaces of stables vector bundles over abelian surfaces}, Nagoya Math. J. 71 1980 
               47-60.                  

\end{thebibliography}
\end{document}